\documentclass[twoside]{article}
\usepackage{amssymb}
\usepackage{graphics}

\oddsidemargin=\evensidemargin
\newtheorem{proposition}{Proposition}
\newtheorem{remark}{Remark}
\newtheorem{conjecture}{Conjecture}

\catcode`\@=11
\long\def\@makefntext#1{
\protect\noindent \hbox to 3.2pt {\hskip-.9pt  
$^{{\eightrm\@thefnmark}}$\hfil}#1\hfill}		

\def\ps@myheadings{\let\@mkboth\@gobbletwo		
\def\@oddhead{\hbox{}
\rightmark\hfil\eightrm\thepage}   
\def\@oddfoot{}\def\@evenhead{\eightrm\thepage\hfil
\leftmark\hbox{}}\def\@evenfoot{}
\def\sectionmark##1{}\def\subsectionmark##1{}}

\catcode`@=11		      		     
\def\ps@plain{\let\@mkboth\@gobbletwo
     \def\@oddhead{}\def\@oddfoot{\eightrm\hfil\thepage
     \hfil}\def\@evenhead{}\let\@evenfoot\@oddfoot}

\newcounter{sectionc}\newcounter{subsectionc}\newcounter{subsubsectionc}
\renewcommand{\section}[1] {\vspace{12pt}\addtocounter{sectionc}{1} 
\setcounter{subsectionc}{0}\setcounter{subsubsectionc}{0}\noindent 
	{\tenbf\thesectionc. #1}\par\vspace{5pt}}
\renewcommand{\subsection}[1] {\vspace{12pt}\addtocounter{subsectionc}{1} 
	\setcounter{subsubsectionc}{0}\noindent 
	{\bf\thesectionc.\thesubsectionc. 
	{\kern1pt \bfit #1}}\par\vspace{5pt}}
\renewcommand{\subsubsection}[1] {\vspace{12pt}
	\addtocounter{subsubsectionc}{1}
	\noindent
	{\tenrm\thesectionc.\thesubsectionc.\thesubsubsectionc.	{\kern1pt 
	\it #1}}\par\vspace{5pt}}
\newcommand{\nonumsection}[1] {\vspace{12pt}\noindent{\tenbf #1}
	\par\vspace{5pt}}

\newcounter{appendixc}
\newcounter{subappendixc}[appendixc]
\newcounter{subsubappendixc}[subappendixc]

\renewcommand{\appendix}[1] {\vspace{12pt}	
	\refstepcounter{appendixc}		
	\setcounter{figure}{0}
	\setcounter{table}{0}
	\setcounter{lemma}{0}
	\setcounter{theorem}{0}
	\setcounter{corollary}{0}
	\setcounter{definition}{0}
	\setcounter{equation}{0}
	\renewcommand{\thefigure}{\Alph{appendixc}.\arabic{figure}}
	\renewcommand{\thetable}{\Alph{appendixc}.\arabic{table}}
	\renewcommand{\theappendixc}{\Alph{appendixc}}
	\renewcommand{\thelemma}{\Alph{appendixc}.\arabic{lemma}}
	\renewcommand{\thetheorem}{\Alph{appendixc}.\arabic{theorem}}
	\renewcommand{\thedefinition}{\Alph{appendixc}.\arabic{definition}}
	\renewcommand{\thecorollary}{\Alph{appendixc}.\arabic{corollary}}
	\renewcommand{\theequation}{\Alph{appendixc}.\arabic{equation}}
	\noindent{\tenbf Appendix \theappendixc #1}\par\vspace{5pt}}

\newcommand{\textlineskip}{\baselineskip=13pt}
\newcommand{\smalllineskip}{\baselineskip=10pt}

\newcommand{\copyrightheading}[1]
	{\vspace*{-2.5cm}\smalllineskip{\flushleft
	{\footnotesize Journal of Knot Theory and Its Ramifications #1}\\
   	{\footnotesize \copyright\kern2pt World Scientific 
         Publishing Company}\\
         }}

\newcommand{\publisher}[2]{{\begin{center}\footnotesize\smalllineskip 
	Received #1\\
	Revised #2
        \end{center}
	}}

\def\abstracts#1#2#3#4{{
	\centering{\begin{minipage}{4.5in}\footnotesize\baselineskip=10pt
	\centerline{ABSTRACT} 
	\parindent=15pt #1\par 
	\parindent=15pt #2\par
	\parindent=15pt #3\par
	\parindent=15pt #4\par
	\end{minipage}}\par}} 

\def\keywords#1{{ 
	\centering{\begin{minipage}{4.5in}\footnotesize\baselineskip=10pt
	{\footnotesize\it Keywords}\/: #1
	\end{minipage}}\par}}

\renewenvironment{thebibliography}[1]
	{\frenchspacing
	 \ninerm\baselineskip=11pt
	 \begin{list}{[\arabic{enumi}]}
	{\usecounter{enumi}\setlength{\parsep}{0pt}
	 \setlength{\leftmargin 13.7pt}{\rightmargin 0pt} 
	 \setlength{\itemsep}{0pt} \settowidth
	{\labelwidth}{[#1]}\sloppy}}{\end{list}}

\newcounter{itemlistc}
\newcounter{romanlistc}
\newcounter{alphlistc}
\newcounter{arabiclistc}

\newcommand{\fcaption}[1]{
        \refstepcounter{figure}
        \setbox\@tempboxa = \hbox{\footnotesize Fig.~\thefigure. #1}
        \ifdim \wd\@tempboxa > 5in
           {\begin{center}
        \parbox{5in}{\footnotesize\smalllineskip Fig.~\thefigure. #1}
            \end{center}}
        \else
             {\begin{center}
             {\footnotesize Fig.~\thefigure. #1}
              \end{center}}
        \fi}

\newcommand{\tcaption}[1]{
        \refstepcounter{table}
        \setbox\@tempboxa = \hbox{\footnotesize Table~\thetable. #1}
        \ifdim \wd\@tempboxa > 5in
           {\begin{center}
        \parbox{5in}{\footnotesize\smalllineskip Table~\thetable. #1}
            \end{center}}
        \else
             {\begin{center}
             {\footnotesize Table~\thetable. #1}
              \end{center}}
        \fi}

\def\pmb#1{\setbox0=\hbox{#1}
	\kern-.025em\copy0\kern-\wd0
	\kern.05em\copy0\kern-\wd0
	\kern-.025em\raise.0433em\box0}

\def\fnm#1{$^{\mbox{\scriptsize #1}}$}		
\def\fnt#1#2{\footnotetext{\kern-.3em
	{$^{\mbox{\scriptsize #1}}$}{#2}}}

\def\fpage#1{\begingroup
\voffset=.3in
\thispagestyle{empty}\begin{table}[b]\centerline{\footnotesize #1}
	\end{table}\endgroup}

\def\runninghead#1#2{\pagestyle{myheadings}
\markboth{{\protect\footnotesize\it{\quad #1}}\hfill}
{\hfill{\protect\footnotesize\it{#2\quad}}}}

\font\tenrm=cmr10
 
\font\tenbf=cmbx10
\font\bfit=cmbxti10 at 10pt
\font\ninerm=cmr9

\font\eightrm=cmr8

\newtheorem{theorem}{Theorem}   

\newtheorem{lemma}{Lemma}

\newtheorem{definition}{Definition}
\newtheorem{corollary}{Corollary}
\def\@begintheorem#1#2{\trivlist	
	\item[\hskip\labelsep{\bf #1\ #2.}]} 
\def\@opargbegintheorem#1#2#3{\trivlist
	\item[\hskip\labelsep{\bf #1\ #2\ (#3).}]}

\newenvironment{proof}{\begin{trivlist}
	\item[\noindent]{\it Proof.}}{\quad $\square$\end{trivlist}} 

    	{\setcounter{itemlistc}{0}		
	 \begin{list}{$\bullet$}		
	{\usecounter{itemlistc}			
	 \leftmargin10pt	       
	 \setlength{\parsep}{0pt}
	 \setlength{\itemsep}{0pt}     
	}}{\end{list}}

\newenvironment{romanlist2}[1]			
	{\setcounter{romanlistc}{0}		
	 \begin{list}{$($\roman{romanlistc}$)$}	
	{\usecounter{romanlistc}		
	 \leftmargin18pt 
	 \setlength{\parsep}{0pt}
	 \setlength{\itemsep}{0pt}	
	 \settowidth{\labelwidth}{#1}                          
	}}{\end{list}}

	{\setcounter{enumii}{0}			
	 \begin{list}{$($\alph{enumii}$)$}	
	{\usecounter{enumii}			
	 \leftmargin18pt		
	 \setlength{\parsep}{0pt}
	 \setlength{\itemsep}{0pt}	
	 \settowidth{\labelwidth}{#1}                          
	}}{\end{list}}

\def\qed{\hbox{${\vcenter{\vbox{			
   \hrule height 0.4pt\hbox{\vrule width 0.4pt height 6pt
   \kern5pt\vrule width 0.4pt}\hrule height 0.4pt}}}$}}

\def\theequation{\thesectionc.\arabic{equation}}  

\begin{document}
\setlength{\textheight}{7.7truein}  

\runninghead{S. Nelson }
{Virtual Crossing Realization}

\normalsize\textlineskip
\thispagestyle{empty}
\setcounter{page}{1}

\copyrightheading{}		    

\vspace*{0.88truein}

\fpage{1}
\centerline{\bf VIRTUAL CROSSING REALIZATION}
\baselineskip=13pt

\centerline{\footnotesize SAM NELSON\fnm{1}\fnt{1}{Permanent address: P.O. Box 2311, Rancho Cucamonga, CA 91729-2311}}
\baselineskip=12pt
\centerline{\footnotesize\it Department of Mathematics}
\baselineskip=10pt
\centerline{\footnotesize\it Whittier College (visiting)}
\centerline{\footnotesize\it 13406 Philadelphia, P.O. Box 634 }
\centerline{\footnotesize\it Whittier, CA 90608-0634}
\centerline{\footnotesize\it knots@esotericka.org}


\vspace*{0.225truein}
\publisher{}

\vspace*{0.21truein} 
\abstracts{
We study virtual isotopy sequences with classical initial and final diagrams,
asking when such a sequence can be changed into a classical isotopy sequence 
by replacing virtual crossings with classical crossings. An example of a 
sequence for which no such virtual crossing realization exists is given. 
A conjecture on conditions for realizability of virtual isotopy sequences 
is proposed, and a sufficient condition for realizability is found. The 
conjecture is reformulated in terms of 2-knots and knots in thickened surfaces.
}{}{}{}

\vspace*{10pt}

\vspace*{1pt}\textlineskip	
\section{Introduction}	
\vspace*{-0.5pt}

In \cite{pgv}, it is observed that classical knot theory embeds in 
virtual knot theory, in the sense that if two classical knots are virtually 
isotopic, then they are classically isotopic; this follows from the fact 
that the fundamental quandle (or alternatively the group system) is 
preserved by virtual moves. Since the fundamental quandle is a complete 
invariant for classical knots, any two classical knot diagrams 
related by a sequence of virtual moves must have isomorphic quandles, and
hence must be isotopic in the classical sense. This also follows from theorem 
1 of \cite{ku}, which says that every stable equivalence class of knots in 
thickened surfaces (which are equivalent to virtual knots) has a unique 
irreducible representative. In \cite{kau}, it is suggested that a purely
combinatorial proof for this fact may be instructive.

A na\"\i ve attempt at a constructive proof in terms of Gauss diagrams
initially looks promising; virtual crossings (and hence virtual moves)
do not appear in Gauss diagrams, which include only classical crossings.
Thus, given a virtual isotopy sequence which begins and ends with 
classical diagrams, we may attempt to construct a classical isotopy 
sequence by simply translating the Gauss diagram sequence to knot diagrams.

This strategy fails because unlike the classical crossing-introducing type 
II move, the Gauss diagram type II move does not require the strands being 
crossed to be adjacent in the plane. Further, the Gauss diagram II move 
permits both the direct and reverse II moves on any pair of strands, 
regardless of the orientation of the strands, while at most one of these 
is realizable using only classical diagrams.

\begin{figure}[!ht]
\[ 
\includegraphics[0.9in,0.9in]{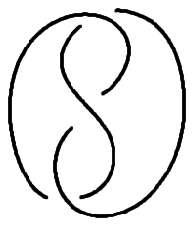} 
\raisebox{0.5in}{$\begin{array}{c}  vII \\ \iff \end{array}$} 
\includegraphics[0.9in,0.9in]{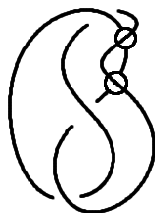} 
\raisebox{0.5in}{$\begin{array}{c}  vII \\ \iff \end{array}$} 
\includegraphics[0.9in,0.9in]{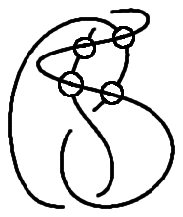} 
\raisebox{0.5in}{$\begin{array}{c}  II \\ \iff \end{array}$} 
\] \[ 
\includegraphics[0.9in,0.9in]{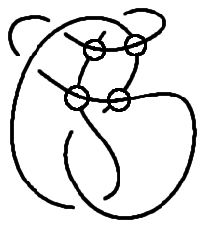} 
\raisebox{0.5in}{$\begin{array}{c}  v \\ \iff \end{array}$} 
\includegraphics[0.9in,0.9in]{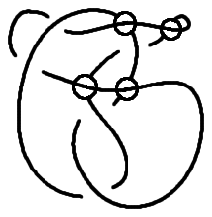} 
\raisebox{0.5in}{$\begin{array}{c}  vI \\ \iff \end{array}$} 
\includegraphics[0.9in,0.9in]{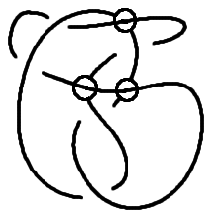} 
\raisebox{0.5in}{$\begin{array}{c}  v \\ \iff \end{array}$} 
\] \[ 
\includegraphics[0.9in,0.9in]{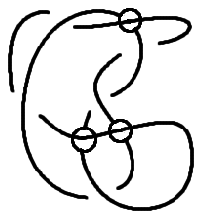} 
\raisebox{0.5in}{$\begin{array}{c}  v \\ \iff \end{array}$} 
\includegraphics[0.9in,0.9in]{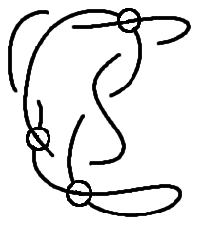} 
\raisebox{0.5in}{$\begin{array}{c}  vI \\ \iff \end{array}$}  
\includegraphics[0.9in,0.9in]{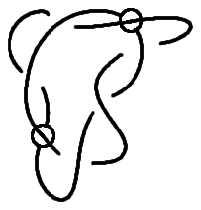} 
\] \[ 
\raisebox{0.5in}{$\begin{array}{c}  II \\ \iff \end{array}$} 
\includegraphics[0.9in,0.9in]{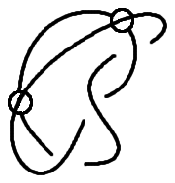} 
\raisebox{0.5in}{$\begin{array}{c}  vII \\ \iff \end{array}$} 
\includegraphics[0.9in,0.9in]{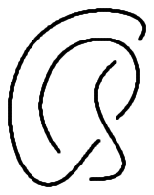} 
\]
\caption{A realized virtual isotopy sequence.}
\label{vcrsequence}
\end{figure}

All of these moves can be realized by introducing virtual crossings. 
Moreover, as these virtual crossings are removed by the end of the sequence, 
it is natural to ask under which circumstances these virtual crossings can 
be replaced with classical crossings throughout the sequence to yield a 
valid classical isotopy sequence. 

An assignment of classical crossing type to each virtual crossing in a virtual 
isotopy sequence is a \textit{virtual crossing realization}. A virtual crossing
realization is \textit{valid} if each classical knot diagram in the resulting 
move sequence differs from the previous diagram by a valid Reidemeister move, 
i.e., if every vII move becomes a II move and every v or vIII move becomes a 
III move. Figure \ref{vcrsequence} depicts a valid virtual crossing 
realization. Here we denote realized virtual crossings as circled classical 
crossings -- these are ordinary classical crossings, the circle is retained 
only to indicate which crossings have been realized. 

Initially we may hope that every virtual isotopy sequence admits a valid 
virtual crossing realization; a combinatorial proof of the type suggested 
in \cite{kau} would then follow. However, as figure \ref{counterexample}
shows, this is not the case. Invariance of the fundamental quandle with
respect to virtual moves implies only the existence of some classical
isotopy sequence between classical diagrams with isomorphic quandles; 
such a sequence might be very different from any given virtual sequence.
In particular, for any pair of equivalent classical diagrams, there 
need only be one classical isotopy sequence with the given end diagrams 
to satisfy the theorem of \cite{pgv}, while virtual sequences with the same 
end diagrams are clearly not unique.

\medskip

In this paper, we study the problem of when a virtual isotopy sequence
may be realized. The paper is organized as follows: We begin with a 
definition of realizability for virtual isotopy sequences. We then 
identify the ways in which a sequence can fail to be realizable, which 
consist of two types of bad moves. We consider each of these types of bad 
moves in turn, proving that one type may always be avoided and analyzing
circumstances in which the other arises. We then obtain our main result, 
theorem \ref{vd}, which gives a sufficient condition for realizability of a
virtual isotopy sequence. A conjecture is proposed, and we reformulate the
conjecture in terms of 2-knots. Finally, the conjecture is reformulated 
in terms of knots in thickened surfaces.

\section{Virtual Knots and Gauss Diagrams}

\begin{figure}[!ht] \label{reidemeister}
\[
\includegraphics[1in,1in]{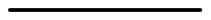}  
\raisebox{0.5in}{$\begin{array}{c} I \\ \iff \end{array}$}
\includegraphics[1in,1in]{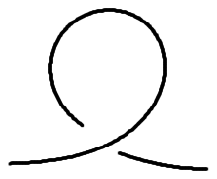}
 \ 
\includegraphics[1in,1in]{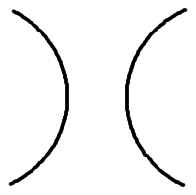}
\raisebox{0.5in}{$\begin{array}{c} II \\ \iff \end{array}$}
\includegraphics[1in,1in]{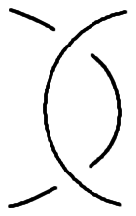}
\]
\[
\includegraphics[1in,1in]{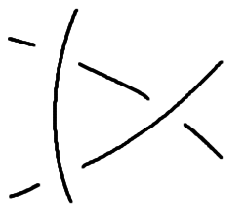}
\raisebox{0.5in}{$\begin{array}{c} III \\ \iff \end{array}$}
\includegraphics[1in,1in]{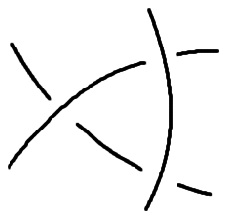}
\]
\caption{Reidemeister moves.}
\label{rmoves}
\end{figure}

A \textit{link diagram} is a planar oriented 4-valent graph with vertices 
regarded as crossings and enhanced with crossing information. The edges are 
oriented so that each vertex has two incoming edges, one over and one under, 
and two outgoing edges, also one over and one under. We may regard knots and 
links combinatorially as equivalence classes of knot and link diagrams under 
the equivalence relation generated by the three Reidemeister moves, pictured 
in figure \ref{rmoves}.

By enlarging the set of decorated graphs to include non-planar 4-valent 
graphs with vertices enhanced with crossing information and edges oriented 
as before, we obtain \textit{virtual links} as equivalence classes under 
the equivalence relation generated by the three Reidemeister moves.

To draw non-planar graphs on planar paper, we must introduce \textit{virtual 
crossings,} which we distinguish from the decorated vertices (or 
\textit{classical crossings}) by denoting virtual crossings as circled 
intersections.

Since these virtual crossings are artifacts of representing non-planar 
graphs in the plane, we may obtain an equivalent diagram by replacing any 
arc containing only virtual crossings with any other arc containing only 
virtual crossings with the same endpoints. This breaks down into 
four virtual moves, one move for each type of thing we can move the arc 
past: the arc itself (move vI), another arc (move vII), a virtual crossing 
(move vIII) and a classical crossing (move v). We may then consider virtual 
knots as equivalence classes of virtual knot diagrams under the equivalence 
relation generated by moves I, II, III, vI, vII, vIII and v, known as 
\textit{virtual isotopy}.

\begin{figure}[!ht]
\[
\includegraphics[0.9in,0.9in]{i1.eps}  
\raisebox{0.5in}{$\begin{array}{c} vI \\ \iff \end{array}$}
\includegraphics[0.9in,0.9in]{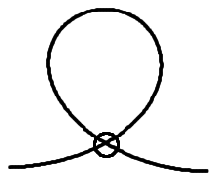}
 \ 
\includegraphics[0.9in,0.9in]{ii1.eps}
\raisebox{0.5in}{$\begin{array}{c} vII \\ \iff \end{array}$}
\includegraphics[0.9in,0.9in]{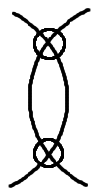}
\]
\[
\includegraphics[0.9in,0.9in]{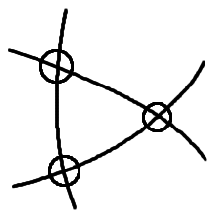}
\raisebox{0.5in}{$\begin{array}{c} vIII \\ \iff \end{array}$}
\includegraphics[0.9in,0.9in]{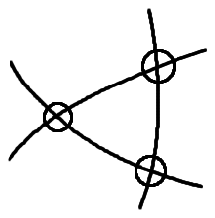}
 \ 
\includegraphics[0.9in,0.9in]{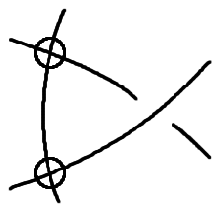}
\raisebox{0.5in}{$\begin{array}{c} v \\ \iff \end{array}$}
\includegraphics[0.9in,0.9in]{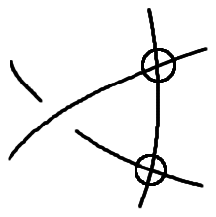}
\]
\caption{Virtual moves.}
\label{vmoves}
\end{figure}

\begin{figure}[!ht]
\[
\includegraphics[0.9in,0.9in]{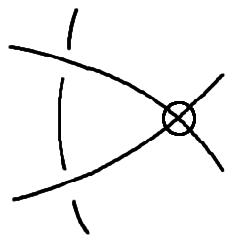}  
\raisebox{0.5in}{$\begin{array}{c} F_h \\ \iff \end{array}$}
\includegraphics[0.9in,0.9in]{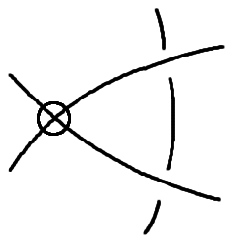}
 \ \ 
\includegraphics[0.9in,0.9in]{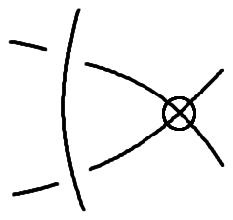}
\raisebox{0.5in}{$\begin{array}{c} F_t \\ \iff \end{array}$}
\includegraphics[0.9in,0.9in]{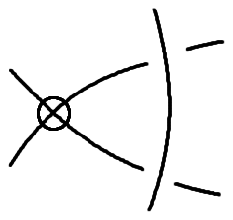}
\]
\caption{Forbidden moves.}
\label{forbidden}
\end{figure}

Two potential moves are not allowed, the two \textit{forbidden moves} 
$F_t$ and $F_h$ (depicted in figure \ref{forbidden}), variants of the 
type III move with two classical crossings and one virtual crossing. 
Unlike the valid virtual moves, the forbidden moves alter the underlying 
graph of the diagram. Together, the two forbidden moves can be used to 
unknot any knot, virtual or classical.

\textit{Gauss diagrams} provide another way of representing virtual knots 
combinatorially. A Gauss diagram for a knot is a circle with oriented 
chords representing crossings; if we think of the circle as the preimage
of the knot diagram under the embedding into $\mathbb{R}^3$ and projection 
to the plane, then the chords join the two preimages of each crossing point. 
We orient the chords ``in the direction of gravity,'' that is, toward the 
preimage of the undercrossing, and we decorate these arrows with signs
given by the local writhe number of the crossing. The Gauss diagram 
of a link has one circle for each component of the link, and crossings 
between components correspond to arrows joining the circles.

\begin{figure}[!ht] 
\[
\includegraphics[0.9in,0.9in]{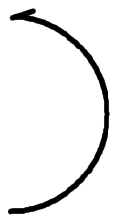}  
\raisebox{0.5in}{$\begin{array}{c} I \\ \iff \end{array}$}
\includegraphics[0.9in,0.9in]{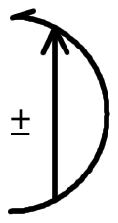}
 \ 
\includegraphics[0.9in,0.9in]{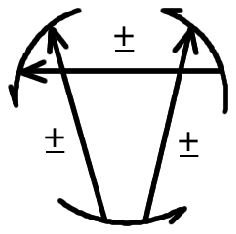}
\raisebox{0.5in}{$\begin{array}{c} III \\ \iff \end{array}$}
\includegraphics[0.9in,0.9in]{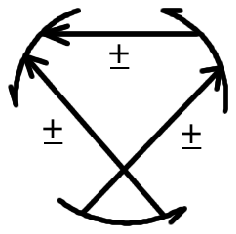}
\]
\[
\includegraphics[1in,1in]{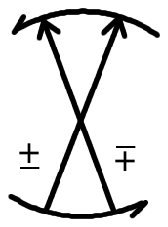}
\raisebox{0.5in}{$\begin{array}{c} II \\ \iff \end{array}$}
\includegraphics[1in,1in]{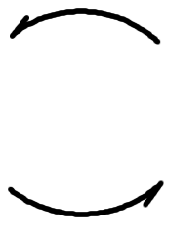}
\raisebox{0.5in}{$\begin{array}{c} II \\ \iff \end{array}$}
\includegraphics[1in,1in]{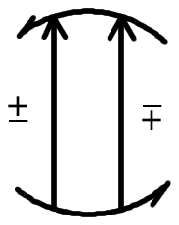}
\]
\caption{Gauss diagram moves.}
\label{gauss}
\end{figure}

Virtual knots and links may then be regarded as equivalence classes of Gauss 
diagrams under the Gauss diagram versions of the Reidemeister moves;  
a Gauss diagram determines a virtual knot diagram up to virtual moves 
(vI, vII, vIII and v), while a virtual knot diagram determines a unique
Gauss diagram. There are several instances of type III moves depending on
the orientation and cyclic order of the three strands involved; only two of 
these are listed in figure \ref{gauss}.

\section{Virtual Isotopy Sequences}

\begin{definition}
A sequence of virtual knot diagrams $K_1 \to K_2 \to \dots \to K_n$ where
$K_i$ differs from $K_{i-1}$ by a single virtual move (and possibly a planar
isotopy) is a \textit{virtual isotopy sequence}. We will sometimes use the
term \textit{valid virtual isotopy sequence} to distinguish a virtual isotopy
sequence from a sequence of virtual knot diagrams in which one or more pairs
of diagrams $K_i$ and $K_{i-1}$ are not related by virtual moves; such a 
sequence may be called an \textit{invalid sequence}.
\end{definition}

The Gauss diagram type II move permits, on any two sections of the circle, 
both the direct type II move, in which both strands are oriented in the same 
direction, and the reverse type II move, in which the strands are oriented 
in opposite directions. Unlike the classical type II move, the Gauss diagram 
II move does not require the arcs being crossed to be adjacent in the 
plane. For any pair of strands, of the four possible Gauss diagram type II 
moves, at most two are classically realizable, and then only if the strands 
are adjacent. However, all of these non-classical moves are realizable in 
virtual knot diagrams with the addition of virtual crossings.

\begin{figure}[!ht] 
\[ 
\includegraphics[0.9in,0.9in]{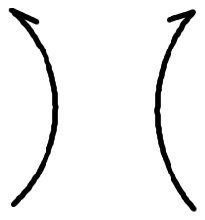}
\raisebox{0.5in}{$\begin{array}{c}  \mathrm{vI} \\ \iff \end{array}$}
\includegraphics[0.9in,0.9in]{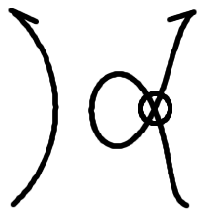} 
\raisebox{0.5in}{$\begin{array}{c} \mathrm{II} \\ \iff \end{array}$} 
\includegraphics[0.9in,0.9in]{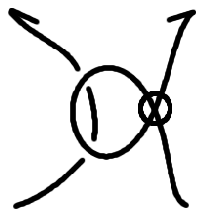} 
\]
\[ 
\includegraphics[0.9in,0.9in]{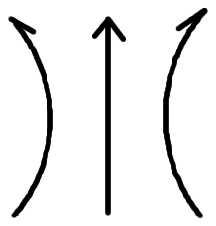} 
\raisebox{0.5in}{$\begin{array}{c} \mathrm{vII} \\ \iff \end{array}$}
\includegraphics[0.9in,0.9in]{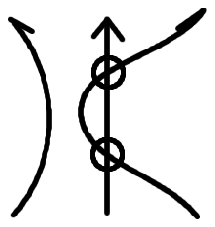}  
\raisebox{0.5in}{$ \begin{array}{c} \mathrm{II} \\ \iff \end{array}$}
\includegraphics[0.9in,0.9in]{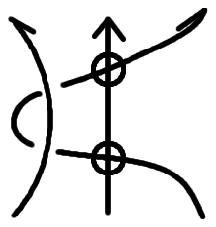} 
\]
\caption{Examples of classically unrealizable type II moves.}
\label{unreal}
\end{figure}

Though a Gauss diagram sequence beginning and ending with realizable classical 
diagrams may not be classically realizable, each of the individual 
unrealizable moves is realizable as a sequence of classical moves on 
realizable classical diagrams if we introduce classical crossings in the 
moves pictured in figure \ref{unreal} instead of virtual crossings. 

Given a virtual isotopy sequence, we can name and follow each virtual 
crossing through the sequence. Introducing classical crossings in place of
virtual crossings then means assigning classical crossing 
information to each virtual crossing throughout the sequence. For each 
individual move, it is clear that we can obtain a legitimate move in this 
way. However, a choice of classical crossing type which makes one move work
may render another move invalid later in the sequence, since crossings 
cannot change type once introduced.

\begin{definition}
An assignment of a sign (+ or -) to each virtual crossing throughout the
diagram is a \textit{virtual crossing realization.} A virtual crossing 
realization is \textit{valid} if it results in a sequence of valid classical 
moves. 
\end{definition}

There are two ways a virtual crossing realization can yield invalid 
moves. One of these is the three-crossing move with all three edges 
alternating, known as the $\Delta$ move; the other is the two-crossing 
move with both edges alternating, called a $\Gamma$ move by analogy with 
the $\Delta$ move.\fnm{2}\fnt{2}{The $\Gamma$ move is sometimes called a 
``2-move.''} Both the $\Delta$ move and $\Gamma$ move are invalid, meaning 
they are not realizable as ambient isotopies. The effect of the $\Delta$ 
move is studied in \cite{D}, and a $\Gamma$ move combined with a pair of 
type II moves effects a crossing change, resulting in unknotting.

\begin{figure}[!ht]
\[ 
\includegraphics[0.9in,0.9in]{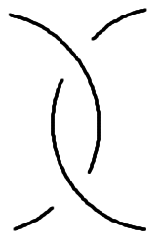}
\raisebox{0.5in}{$ \begin{array}{c}  \Gamma \\ \iff \end{array}$} 
\includegraphics[0.9in,0.9in]{ii1.eps} 
\]
\[ 
\includegraphics[0.9in,0.9in]{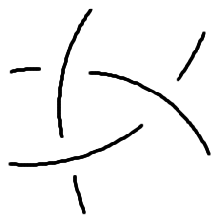} 
\raisebox{0.5in}{$\begin{array}{c}  \Delta \\ \iff \end{array}$} 
\includegraphics[0.9in,0.9in]{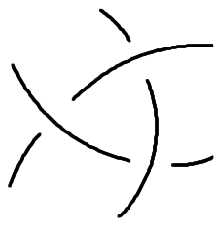} 
\]
\caption{Invalid moves arising in virtual crossing realization.}
\label{invalid}
\end{figure}

Thus, a virtual crossing realization with no $\Delta$ or $\Gamma$ moves 
yields a classical isotopy sequence from the initial diagram to the final 
diagram. Thus we ask, which virtual isotopy sequences admit a valid virtual 
crossing realization? 

\begin{theorem}
Not every valid virtual isotopy sequence with classical end diagrams admits
a valid virtual crossing realization. \label{countertheorem}
\end{theorem}

\begin{figure}[!ht]
\[
\includegraphics[1in,1in]{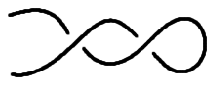}
\raisebox{0.5in}{$\begin{array}{c} vI \\ \iff \end{array}$} 
\includegraphics[1in,1in]{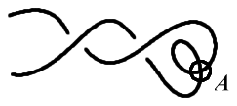} 
\raisebox{0.5in}{$\begin{array}{c} vII \\ \iff \end{array}$} 
\includegraphics[1in,1in]{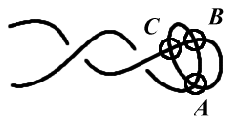}  
\]\[
\raisebox{0.5in}{$\begin{array}{c}  v \\ \iff \end{array}$} 
\includegraphics[1in,1in]{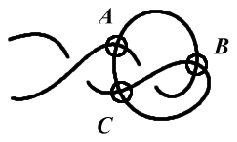} 
\raisebox{0.5in}{$\begin{array}{c}  v \\ \iff \end{array}$} 
\includegraphics[1in,1in]{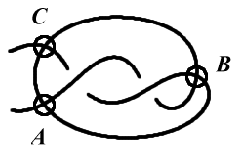} 
\raisebox{0.5in}{$\begin{array}{c}  v \\ \iff \end{array}$} 
\includegraphics[1in,1in]{c4.eps} 
\]\[
\raisebox{0.5in}{$\begin{array}{c}  v \\ \iff \end{array}$} 
\includegraphics[1in,1in]{c3.eps} 
\raisebox{0.5in}{$\begin{array}{c}  vII \\ \iff \end{array}$} 
\includegraphics[1in,1in]{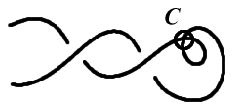} 
\raisebox{0.5in}{$\begin{array}{c}  vI \\ \iff \end{array}$} 
\includegraphics[1in,1in]{c1.eps} 
\]
\caption{A virtual isotopy sequence which does not admit a virtual crossing realization.}
\label{counterexample}
\end{figure}

\begin{proof}
Figure \ref{counterexample} depicts 
a virtual isotopy sequence with no valid virtual crossing realization; 
each of the eight possible assignments of signs to the virtual crossings 
yields either a $\Delta$ or a $\Gamma$ move. Moreover, this sequence
may be spliced in to any other virtual isotopy sequence to yield a new
sequence which does not admit a virtual crossing realization.
\end{proof}

However, changing either of the classical crossings in the sequence in 
figure \ref{counterexample} yields a valid sequence which does admit a 
valid virtual crossing realization, and whose end diagrams are equivalent 
to the originals. Indeed, the sequence of figure \ref{counterexample} may 
be viewed as a truncation of a longer sequence which has been partially
realized; in this scenario, the two classical crossings have been 
introduced, and are later removed, in realized type vI moves. Moreover,
if we choose the opposite sign for either of these classical 
crossings, the resulting sequence admits a realization.

Indeed, let $S$ be a virtual isotopy sequence which admits a virtual crossing 
realization. Then we may obtain a virtual isotopy sequence which 
does not admit a virtual crossing realization by realizing some crossings 
in $S$ and then truncating the sequence provided the realizations chosen for 
the realized crossings yield a valid virtual isotopy sequence yet are 
incompatible with the remaining possible choices, as in figure 
\ref{counterexample}.

Thus, in order to find a classical isotopy sequence given a virtual sequence 
with classical end diagrams, we must avoid situations like the one in
figure \ref{counterexample}. We 
begin by considering which classical crossings may be replaced with virtual 
crossings, yielding an equivalent sequence.

\section{Classical Crossing Virtualization}

\begin{definition}
A classical crossing in a virtual knot diagram is \textit{virtualizable} if 
the virtual knot diagram obtained by replacing the crossing with a virtual 
crossing is virtually isotopic to the original diagram. More generally, a 
set of crossings in a virtual knot diagram is virtualizable if the diagram
obtained by replacing each crossing in the set with a virtual crossing is
equivalent to the original diagram. Similarly, a set of classical crossings
is \textit{switchable} if the knot diagram obtained by switching the crossing 
type of all crossings in the set is equivalent to the original diagram. 
\end{definition}

\begin{figure}[!ht]
\[
\includegraphics[1in,1in]{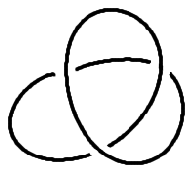} 
\raisebox{0.5in}{$\iff$} 
\includegraphics[1in,1in]{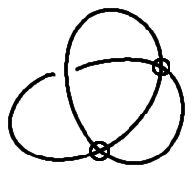} \hskip 0.2in
\raisebox{0.5in}{$\not\hskip -3mm\iff$} 
\includegraphics[1in,1in]{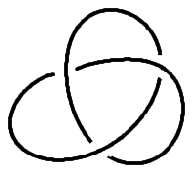} 
\]
\caption{Virtualizable crossings need not be switchable.}
\label{switch}
\end{figure}

Switchability does not imply virtualizability; if we simultaneously switch 
all the crossings in a square knot, for example, we obtain an equivalent 
diagram, but virtualizing all the crossings yields an unknot. Similarly, 
virtualizability does not imply switchability. Figure \ref{switch} shows a 
virtual knot digram with a pair of crossings which are virtualizable 
but not switchable.

A crossing which appears in neither the initial nor the final diagram in 
a virtual isotopy sequence is \textit{temporary}. Temporary crossings are 
both introduced and removed during the course of the isotopy, and may be 
classical or virtual. In a virtual isotopy sequence whose initial and final 
diagrams contain only classical crossings, every virtual crossing is 
temporary; such a sequence may also contain classical temporary 
crossings.

Virtualizability and switchability are properties of sets of crossings, and
may be defined both for individual diagrams and for sequences of diagrams.
Intuitively, a set of crossings is virtualizable (respectively, switchable)
in a diagram if virtualizing (resp., switching) the crossings in the diagram
results in an equivalent diagram. Similarly, a set of crossings in 
virtualizable (respectively, switchable) in a virtual isotopy sequence if 
virtualizing (resp., switching) the crossings in the diagram results in an 
equivalent sequence. More formally, we have the following:

\begin{definition}
Let $K=K_1\to\dots\to K_n$ be a virtual isotopy sequence and let $J$ be a
set of crossings in $K$. Then the set $J$ is \textit{sequentially 
virtualizable} if the move sequence $K'$ obtained by replacing each crossing 
in $J$ throughout $K$ with a virtual crossing is a valid virtual isotopy 
sequence, with $K_1'$ equivalent to $K_1$ and $K_n'$ equivalent to $K_n$. 
\end{definition}

\begin{proposition}
A set of crossings $J$ in a virtual isotopy sequence $K$ is sequentially 
virtualizable if and only if $J$ satisfies the following conditions: 
\begin{romanlist2}{(ii)}
\item The subset of non-temporary crossings in $J$ is 
virtualizable in each end diagram, and 
\item No crossing in $J$ appears in a type III move with 
two crossings not in $J$.
\end{romanlist2} 
\end{proposition}

\begin{proof}
Condition (i) says that the sequence $K'$ obtained from $K$ by replacing 
the crossings in $J$ with virtual crossings has end diagrams which are 
equivalent to the end diagrams of $K$. Condition (ii) says that type III
moves in $K$ get replaced with type v or vIII moves in $K'$; in 
particular, $K'$ is free from forbidden moves.

Conversely, if virtualizing the crossings in $J$ results in a valid
virtual isotopy sequence, the absence of forbidden moves implies condition
(ii), and $K_1$ equivalent to $K_1'$ and $K_n$ equivalent to $K_n'$ imply
that any non-temporary crossings are virtualizable.
\end{proof}

Similarly, we have

\begin{definition}
Let $K=K_1\to\dots\to K_n$ be a virtual isotopy sequence and let $J$ be a
set of classical crossings in $K$. Then the set $J$ is \textit{sequentially
switchable} if the move sequence $K'$ obtained by switching each crossing 
in $J$ throughout $K$ is a valid virtual isotopy sequence, with $K_1'$ 
equivalent to $K_1$ and $K_n'$ equivalent to $K_n$. 
\end{definition}

\begin{proposition}
A set of crossings is sequentially switchable iff 
\begin{romanlist2}{(ii)}
\item Every non-temporary crossing in the set is switchable in the end 
diagrams and
\item If one crossing in the set lies on an over-arc in a III move, the 
other crossing on the over-arc is also in the set.
\end{romanlist2} 
\end{proposition}

\begin{proof}
Switching one crossing on an over-arc changes a III move to a $\Delta$ move, 
but switching both keeps it a type III.
\end{proof}

\begin{definition}
A virtual isotopy sequence is \textit{maximally virtualized} if it contains
no sequentially virtualizable classical crossings. A virtual isotopy sequence
\textit{has switch-free ends} if every sequentially switchable set of crossings
is virtualizable.
\end{definition}

A union of sequentially virtualizable sets of crossings is sequentially 
virtualizable; hence, in any virtual isotopy sequence there is a maximal 
sequentially virtualizable set. Given a virtual isotopy sequence $K$, we 
may obtain a maximally virtualized sequence by replacing every crossing in 
the maximal sequentially virtualizable set with a virtual crossing. Note 
that a subsequence of a maximally virtualized sequence need not be maximally 
virtualized.

\begin{remark}
Though our interest in classical crossing virtualization is motivated by
virtual crossing realization, it may occasionally be desirable to reduce the 
number of crossings needed in Gauss diagram isotopy sequences. A maximally 
virtualized isotopy sequence has the minimal number of classical crossings
among virtual isotopy sequences with the specified underlying planar graph 
sequence.
\end{remark}

A virtual isotopy sequence which contains sequentially virtualizable 
classical crossings may be viewed as a partially completed virtual crossing 
realization problem. Such a sequence is ``correctly'' partially completed in 
the sense that no invalid moves have yet been introduced by the choice of 
realization for the virtualizable crossings; however, a virtual isotopy 
sequence like the one in figure \ref{counterexample} with virtualizable 
classical crossings which does not admit a virtual crossing realization 
might admit one if we are allowed to switch the sequentially virtualizable 
classical crossings.

\begin{conjecture}\label{main}
Every maximally virtualized virtual isotopy sequence with switch-free 
classical initial and final diagrams admits a valid virtual crossing 
realization.
\end{conjecture}

Note that a maximally virtualized sequence with classical end diagrams
has no sequentially virtualizable non-temporary crossings. Since sequentially 
virtualizable crossings in a sequence are sequentially virtualizable in every 
subsequence, requiring the end diagrams to be classical and free of 
switchable virtualizable crossings avoids situations such as figure 
\ref{counterexample} in which a realizable sequence is made unrealizable 
by truncation.

\section{$ir$ Classes and Realization Sets}

\begin{definition}
A \textit{realization set} is a set of virtual crossings in a virtual 
isotopy sequence for which a valid virtual crossing realization exists. 
\end{definition}

Conjecture \ref{main} says that if a virtual isotopy sequence is maximally
virtualized and has switch-free classical end diagrams, then the set of all
virtual crossings in the set is a realization set.

Clearly, if a set of virtual crossings is obtained by virtualizing a 
sequentially virtualizable set of classical crossings, it is a realization
set; as figure \ref{counterexample} shows, not every set of virtual
crossings is a realization set. We now consider when a set of virtual 
crossings in a virtual isotopy sequence is a realization set.

Let $X$ be the set of virtual and sequentially virtualizable crossings in a 
virtual isotopy sequence which starts and ends with realizable classical 
diagrams. For each $x\in X$, define $i(x)=x$ if $x$ is either present in the 
initial diagram or is introduced in a type I or vI move; otherwise, $x$ is 
introduced in a type II or vII with another crossing $y$, in which case set 
$i(x)=y$. Similarly, define $r(x)=x$ if $x$ is either present in the final 
diagram or removed in a type I or vI move; otherwise, $x$ is removed in a 
type II or vII with another crossing $y$, in which case set $r(x)=y$. We 
then have involutions $i:X\to X$ and $r:X\to X$ taking each crossing 
$x\in X$ to its \textit{introduction partner} $i(x)$ and each crossing 
$x\in X$ to its \textit{removal partner} $r(x)$. In particular, both $i$ and 
$r$ are injective. For any $x\in X$, we may have $i(x)=r(x)$ or $i(x)\ne 
r(x)$. Reversing the order of steps in the isotopy sequence interchanges 
$i$ and $r$.

The equivalence classes of sequentially virtualizable crossings under the 
equivalence relation generated by the relations $x \sim i(x)$ and $x\sim 
r(x)$ are \textit{$ir$ classes}. The set of $ir$ classes forms a partition 
on the set $X$ of virtual and sequentially virtualizable crossings in a 
virtual isotopy sequence. We can represent an $ir$ class graphically with an 
\textit{$ir$ diagram} as follows: an $ir$ diagram is a graph with a vertex 
for each crossing in the $ir$ class, an edge labeled $i$ joining $x$ to 
$i(x)$ and an edge labeled $r$ joining $x$ and $r(x)$ for each crossing 
$x$ in the $ir$ class.
 
\begin{figure}[!ht]
\[
\includegraphics[2in,1in]{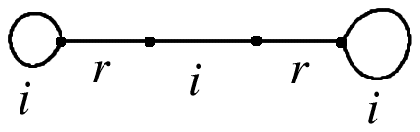} 
\]
\[
\includegraphics[2in,1in]{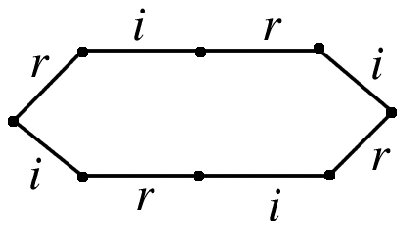}
\]
\caption{$ir$ class diagrams.}
\label{ird}
\end{figure}

\begin{proposition}
Classical crossings may be virtualized to obtain a valid virtual 
isotopy sequence only by virtualizing $ir$ classes. A union of $ir$ classes
of classical crossings is sequentially virtualizable only if no crossing 
in any class in the set appears in a type III move with two classical
crossings not in any class in the set. In particular, a lone $ir$ class is 
sequentially virtualizable only if no crossing in the class appears in a 
type III move with two classical crossings not in the class.
\end{proposition}

\begin{proof}
Virtualizing $x$ but not $i(x)\ne x$ (or $r(x)\ne x$) results in an invalid 
pseudo-II move with one classical and one virtual crossing. If any crossing 
in the set of $ir$ classes appears in a type III move with two classical 
crossings not in the set, virtualizing that crossing will change the  
III move into an invalid move, either one of the two forbidden moves $F_t$ 
or $F_h$ of figure \ref{forbidden} or an invalid move equivalent to one of 
the two forbidden move sequences $F_o$ or $F_s$ in \cite{nel}. 

Virtualizing a complete $ir$ class of temporary classical crossings
changes the type I and II moves to valid vI and vII moves, and the condition
that no crossing in the set of $ir$ classes being virtualized appears in a 
type III move with two crossings not in any class in the set implies that 
each type III move is virtualized to either a valid type v move or a valid 
type vIII move, that is, $\Delta$ moves are avoided.
\end{proof}

A realization set therefore must be a union of $ir$ classes. Moreover,
a realization set must must be ``closed under type v and vIII'' moves, in the
sense that if any crossing in the realization set appears in a type v move, 
the $ir$ class of the other virtual crossing must also be included in the 
set, and no crossing in any class in the set may appear in a type vIII move 
with two crossings whose classes are not in the set, in order to avoid 
forbidden moves.

\begin{proposition} \label{ir}
An $ir$ class has either
\begin{romanlist2}{(iii)}
\item $i(x)=x=r(x)$, or
\item $i(x)=x$, $r(y)=y$ for a unique $x$ and $y$, $x\ne y$, or 
\item $i(x)\ne x$ and $r(x)\ne x$ for all $x$ in the class.
\end{romanlist2}
The $ir$ diagram of (\textit{i}) is a single vertex with two loops, one 
labeled $i$ and one labeled $r$. The diagram of (\textit{ii}) is a sequence 
of vertices connected by single edges with a loop at each end. The diagram 
of (\textit{iii}) is a closed loop with an even number of edges and an 
even number of vertices, with edges alternately labeled $r$ and $i$. 
\end{proposition}

\begin{proof}
The maps $i$ and $r$ are injective, so every vertex meets one $i$ edge and 
one $r$ edge. If $i(x)=x$ or $r(x)=x$ for a vertex $x$, the $ir$ diagram has a 
loop at $x$; the other edge at $x$ can be another loop or it can connect to 
another vertex. This next vertex can either have a loop or connect to another 
new vertex, but it cannot connect back to a previous vertex since each vertex 
already listed in the diagram has met both an $i$ and $r$ edge. After some 
number of steps, we reach the vertex representing the final crossing in the 
$ir$ class, which must therefore meet a loop.

If an $ir$ class has no crossing with $i(x)=x$ or $r(x)=x$, then the 
graph is a closed loop, since every vertex meets one $i$ and one $r$ edge.
Thus the number of vertices in the $ir$ diagram equals the number of edges; 
this number is even since the number of $i$ edges equals the number of $r$ 
edges. Note that the edges must alternate between $r$ and $i$ labels.
\end{proof}

\begin{corollary}
An $ir$ class with an odd number $n\ge 3$ of crossings must have two crossings
either present in the end diagrams or introduced or removed in type I moves.
\end{corollary}

\begin{definition}
An $ir$ \textit{class realization} is a choice of sign for each crossing 
in an $ir$ class. An $ir$ class realization is \textit{valid} if the
resulting sequence does not contain any $\Gamma$ moves.
A \textit{virtual set realization} is a choice of crossing
sign for each virtual crossing in a set.
\end{definition}

\begin{proposition}\label{nogamma}
Each virtual $ir$ class admits two $ir$ class realizations without 
$\Gamma$-moves.
\end{proposition}

\begin{proof}
In a $\Gamma$ move, the crossings have the same sign, while in both the direct 
and reverse II moves, the crossing pairs have opposite signs. Thus, if we 
assign alternating signs to distinct crossings connected by edges in a virtual 
$ir$ class diagram, the resulting virtual crossing realization contains no 
$\Gamma$ moves between crossings in that $ir$ class. For each type of $ir$ 
class diagram, we can make alternating sign assignments consistently. If 
$i(x)=x$ is assigned $\epsilon =\pm 1$, then $r(x)$is assigned $-\epsilon$,
and $i(r(x))$ gets $(-1)^2\epsilon=\epsilon$, etc., until we reach the other 
end of the $ir$ class diagram. In a closed loop $ir$ class diagram with $2k$ 
vertices, choose a starting vertex $x$ and assign it $\epsilon$, then assign 
$-\epsilon$ to $i(x)$, $(-1)^2\epsilon=\epsilon$ to $r(i(x))$ and continue 
around the loop; when we reach $x$ again, it gets assigned 
$(-1)^{2k}\epsilon=\epsilon$, and the assignment is consistent.

A choice of sign for one crossing in an $ir$ class thus determines the signs 
for the whole class; hence for each $ir$ class there are two alternating 
assignments of signs, and thus two $ir$ class realizations which do not 
contain $\Gamma$ moves within the $ir$ class. 
\end{proof}

\begin{corollary}\label{2^n}
A union of $n$ $ir$ classes has $2^n$ virtual set realizations without
$\Gamma$ moves. In particular, if a virtual isotopy sequence has $n$ virtual 
$ir$ classes there are $2^n$ virtual crossing realizations which do 
not contain $\Gamma$ moves.
\end{corollary}

Say that a virtual set realization whose restriction to each $ir$ class is 
one of the two valid $ir$ class realizations \textit{respects} $ir$ 
\textit{classes}. Then any valid virtual crossing realization must respect 
$ir$ classes. For a given set of virtual crossings in a virtual isotopy 
sequence, to determine whether the set is a realization set it suffices to
check only the valid $ir$ class realizations. 

\begin{remark}
We can now see why the virtual isotopy sequence in figure \ref{counterexample}
has no valid virtual crossing realization: it has only one $ir$ class, and 
hence only two of the eight virtual crossing realizations are free of 
$\Gamma$ moves. Inspection reveals that both of these realizations make the 
edge connecting crossing $A$ and crossing $C$ alternating; either choice then 
makes one of the two following v moves a $\Delta$. Note that both classical 
crossings in this sequence are sequentially virtualizable, so this example 
does not contradict conjecture \ref{main}.
\end{remark}

\section{$\Delta$ moves}

In this section we fix a virtual isotopy sequence with a choice of virtual 
crossing realization and consider when the realized sequence includes 
$\Delta$ moves.

\begin{definition}
An edge in a virtual link diagram whose endpoints belong to the 
same $ir$ class is an \textit{intra-class} edge. An edge
joining crossings in distinct $ir$ classes is an \textit{inter-class} edge.
A type v or vIII move is \textit{inter-class} if any of the three
edges joining the crossings in the move is an inter-class edge; otherwise,
the move is \textit{intra-class.} 
\end{definition}
 
If an intra-class edge is alternating for one choice of valid 
$ir$-class realization, it is also alternating for the other choice; if an 
intra-class edge is non-alternating for a choice of valid $ir$ class 
realization, it is non-alternating for the other choice. If a $\Delta$ move 
involves an pair of $i$- or $r$-partners in a virtual crossing realization 
which respects $ir$ classes, the edge connecting the pair in the move cannot 
be one of the edges originally joining the pair, since these are 
non-alternating for both 
class-respecting realizations. However, $ir$ partners may be connected in a 
move by non-original edges, and an intra-class move need not contain 
$ir$ partners, only a pair of crossings from the same $ir$ class. 

If an alternating intra-class edge appears in a type v move, then only one 
of the two $ir$ class realizations makes the move a valid III move; in this 
case, the move \textit{determines} a choice of realization for the class. 
Moreover, one determined $ir$ class may determine another, if a crossing from 
the determined class appears opposite an alternating intra-class edge in a 
type vIII move. Note that switching any crossing in a $\Delta$ move 
makes the move valid, while switching either of the two crossings on the 
over-arc in a type III move changes it to a $\Delta$. 

These observations enable us to classify possible counterexamples 
to conjecture \ref{main} into three types.

A counterexample to conjecture \ref{main} is \textit{type i} if its $\Delta$
moves involve crossings from a single $ir$ class. 
It might have an intra-class type vIII move in which all three edges are 
alternating; then the move is realized as a $\Delta$ for both choices 
of $ir$-class realization. Another example of this type 
would have an $ir$ class with two intra-class type v moves with an alternating 
edge in each move, so that the two moves determine opposite realizations for 
the class. The example in figure \ref{counterexample} would be of this type, 
if the sequence were maximally virtualized.

A counterexample to conjecture \ref{main} is \textit{type ii} if its $\Delta$ 
moves involve distinct determined $ir$ classes 
with incompatible $ir$ class realizations. A pair of crossings from distinct 
determined $ir$ classes might meet in an inter-class v move which is realized 
as a $\Delta$ by the determined $ir$ class realizations, or a vIII move could 
have three crossings from determined classes which realize as a $\Delta$. 
Another example might have two determined $ir$ classes which determine 
opposite realizations for a third $ir$ class.

A counterexample to conjecture \ref{main} is \textit{type iii} if its 
$\Delta$ moves involve a chain of $ir$ classes with inter-class vIII moves 
such that every realization respecting $ir$ classes realizes at least one of 
these vIII moves as a $\Delta$. That is, any attempt to fix the move by 
switching one $ir$ class realization simply creates a new inter-class 
$\Delta$ move. An example of this type might have a ``cycle'' of $ir$ classes 
so that fixing a $\Delta$ move by switching one $ir$ class changes another 
III move to a $\Delta$; fixing this move by switching the $ir$ class of 
another crossing in the move then changes another III move to a $\Delta$, 
and so on, until we eventually break the first $\Delta$ move we fixed. A 
minimal example of this type would have three $ir$ classes, say $A,B,C$, and
inter-class moves between each, such that all eight choices of realization
for the triple $\{A,B,C\}$ realize at least one of the inter-class moves as 
a $\Delta$.

Any of these situations in a maximally virtualized virtual isotopy sequence 
with realizable switch-free classical end diagrams would contradict 
conjecture \ref{main}. Conversely, proving that none of the three types of 
situation listed above can occur would establish conjecture \ref{main}.

\section{Virtually Descending Diagrams}

We now give a sufficient condition for when a virtual isotopy 
sequence with single-component classical end diagrams has a valid virtual 
crossing realization. In this section, $K$ is a virtual knot diagram, 
i.e., a single-component virtual link digram.

\begin{definition}
A virtual crossing realization of $K$ is \textit{virtually descending} 
with respect to a chosen base point and orientation if, starting at the base 
point and following the orientation, we encounter each realized virtual 
crossing first as an overcrossing. Say that a move \textit{fixes} a base point 
if the base point lies outside the part of the diagram pictured in the move.
\end{definition}

\begin{lemma}\label{vdiii}
If the base point is fixed by a type III move and $K$ is virtually 
descending before the move, $K$ is virtually descending after the move.
\end{lemma}

\begin{proof}
Consider a realized vIII move. For $K$ to be virtually descending, the 
strands must be encountered in the order listed in the picture below. 
Inspection shows that the diagrams are virtually descending both before and 
after the move.
\[
\includegraphics[0.9in,0.9in]{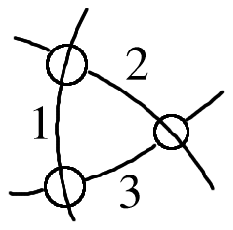} \ \ 
\raisebox{0.5in}{$\begin{array}{c} \mathrm{III} \\ \iff \end{array}$}
\includegraphics[0.9in,0.9in]{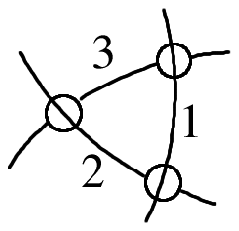} 
\]
The other cases are similar.
\end{proof}

This observation suggests a strategy for choosing virtual crossing 
realizations, namely find a base point fixed by all v, vI, vII and vIII moves, 
then realize the crossings to make the diagram virtually descending. 
Lemma \ref{vdir} shows that this strategy is compatible with our 
previous work.

\begin{lemma} \label{vdir}
If a virtual isotopy sequence fixes a base point, realizing the virtual
crossings to make $K$ virtually descending at each step in the sequence 
respects $ir$ classes.
\end{lemma}

\begin{proof}
Choose an orientation and realize the crossings at each step to keep $K$ 
virtually descending at each move. Since the sequence of moves fixes the base 
point, every edge joining $i$ or $r$ pairs in a vII move is made 
non-alternating; otherwise, the base point lies on one of the edges, contrary 
to assumption. Given such a base point and vII move, the two choices of 
orientation yield the two choices of valid $ir$ class realization from 
proposition \ref{ir}.
\end{proof}

\begin{lemma} \label{vdviii}
If a $K$ is virtually descending before a realized type vIII move with respect 
to a base point fixed by the move, the move is realized as a valid III move.
\end{lemma}

\begin{proof}
If $K$ is virtually descending, then the first strand encountered in the 
move meets both realized crossings going over, hence the edge connecting
them is non-alternating and the move is valid.
\end{proof}

\begin{lemma} \label{vdv}
If $K$ is virtually descending before a type v move with respect to a base 
point fixed by the move, the move is realized as a $\Delta$ move if and
only if the classical undercrossing is encountered before the second virtual 
crossing when following the knot from the chosen base point and orientation. 
That is, the move is a $\Delta$ iff the strand with the classical 
undercrossing is not third in the cyclic ordering of strands determined by
the chosen base point and orientation.
\end{lemma}

\begin{proof}
Consider a realized type v move.
A virtually descending diagram in which both virtual crossings are encountered 
before the classical crossing has both virtual overcrossings adjacent; hence 
the move is not a $\Delta$. A virtually descending diagram in which the 
classical overcrossing is on the first strand encountered likewise has a pair
of adjacent overcrossings, and hence is not a $\Delta$. Thus, for a virtually 
descending diagram to be in position for a $\Delta$ move, the crossings must 
be encountered the order illustrated, namely, classical undercrossing before 
the second virtual crossing.
\[
\includegraphics[0.9in,0.9in]{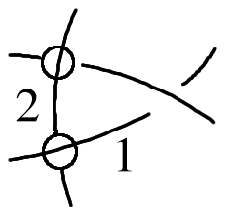}  
\hskip 0.1in \raisebox{0.5in}{$\iff$} 
\includegraphics[0.9in,0.9in]{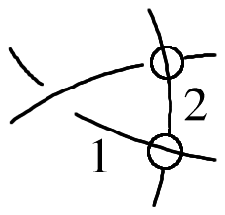} 
\]
Conversely, one checks that the situation illustrated is indeed a $\Delta$ 
move.
\end{proof}

\begin{theorem} \label{vd}
If a virtual isotopy sequence fixes a base point such that the overcrossing 
is always encountered before the undercrossing for all classical
crossings involved in type v moves for a choice of orientation of the knot,
the virtual isotopy sequence admits a valid virtual crossing realization.
\end{theorem}

\begin{proof}
Realize all virtual crossings as they are introduced to make the diagrams 
virtually descending with respect to the given base point and orientation. 
The first possible $\Delta$ move is either of type v or vIII; if the former, 
lemma \ref{vdv} implies the move is realized as a valid type III move,
while if the latter, lemma \ref{vdviii} yields the same conclusion. 
Then lemma \ref{vdiii} implies that the next diagram is virtually
descending; then the next v or vIII move is also a valid III move, and
we repeat until we reach the final diagram. Lemma \ref{vdir} implies
that the sequence is free from $\Gamma$ moves. Thus, the virtual crossing 
realization specified is valid.
\end{proof}

\begin{remark}
The proof of theorem \ref{vd} relies on the fact that the entire diagram 
is virtually descending at every step. Thus, an attempt to apply this method 
to individual $ir$ classes, e.g., considering distinct base points and 
orientations for each $ir$ class, fails in general. If the set of virtual 
crossings in a given virtual isotopy sequence can be divided into disjoint 
non-interacting classes, then distinct base points and orientations satisfying 
the condition of theorem \ref{vd} may yield the result, though such a 
sequence may also be reduced after possible move reordering into shorter 
sequences.
\end{remark}

\section{Knotted Surfaces}

A 2-knot is a compact smooth surface embedded in $\mathbb{R}^4$. A 2-knot 
diagram is a compact smooth surface $M$ immersed in $\mathbb{R}^3$ with 
singular set enhanced with crossing information. As with 1-knots, at a 
crossing curve we indicate which sheet goes over by drawing the undercrossing 
sheet ``broken.''  The 
preimage of the singular set is a set of closed curves and arcs in $M$, 
called the \textit{double decker curves}; double decker arcs end at 
\textit{branch points}. The double decker curves are divided into 
\textit{upper decker curves} on the upper sheet at each crossing and 
\textit{lower decker curves} on the lower sheet, analogous to upper and lower
crossing point preimages in an ordinary 1-knot. Each double point curve is the 
image of an upper decker curve and a lower decker curve. When three sheets 
meet at a triple point, one sheet is highest, one between the others, and one 
lowest.

If we take a 2-dimensional slice of a 2-knot diagram in $\mathbb{R}^3$ by 
intersecting the 2-knot diagram with a plane missing any triple points, 
we obtain an ordinary link diagram; conversely, if we stack the diagrams 
in an isotopy sequence, letting the link diagram sweep out a broken surface 
diagram in $\mathbb{R}^3$, we obtain a portion of a 2-knot diagram connecting 
the initial and final diagrams.\fnm{3}\fnt{3}{More precisely, we have a 
link concordance connecting the link diagrams.} Triple points in the 
resulting 2-knot diagram correspond to Reidemeister III moves. Call the 
direction normal to the planes of the link diagrams \textit{vertical} and the 
planes \textit{horizontal.} Note that taking slices of an arbitrary 2-knot 
diagram does not typically yield an isotopy sequence, since local extrema in 
the vertical direction may result in differing numbers of link components in 
the resulting link diagrams.

Representing virtual crossings as undecorated self-intersections, a virtual 
isotopy sequence sweeps out an immersed broken surface diagram in 
$\mathbb{R}^3$ with singular set divided into classical or ``decorated'' 
(crossing information specified) and virtual or ``undecorated'' (no crossing 
information specified) parts. Conversely, an immersed surface diagram with 
some decorated and some undecorated double point curves represents a virtual 
isotopy sequence only if it has no local extrema in the vertical direction 
and no triple points with one undecorated and two decorated arcs.

An immersed surface which corresponds to a Reidemeister move sequence
can be lifted to a 2-knot diagram, that is, we can choose crossing 
information along the undecorated singular set which makes the immersed 
surface a portion of an ordinary 2-knot diagram. 

\begin{theorem}
Conjecture \ref{main} is equivalent to the following:
If $M$ is an immersed smooth broken surface in $\mathbb{R}^3$ satisfying 
\begin{romanlist2}{(iii)}
\item $M$ has no local extrema in the vertical direction,
\item the boundary of $M$ is a pair of switch-free classical link diagrams, 
\item the crossing information on the sheets of $M$ is compatible with the 
crossing information in the end diagrams and every triple point with three
decorated curves has a highest, middle and lowest sheet,
\item no undecorated double point curve intersects two decorated double point
curves at a triple point, and
\item all double point curves not reaching the end diagrams and meeting at 
most one decorated curve are undecorated,
\end{romanlist2}
then $M$ lifts to a surface knot diagram.
\label{re-1}
\end{theorem}

\begin{proof}
Conditions (i) - (iv) guarantee that $M$ defines a virtual isotopy sequence:
(i) says that the number of components stays constant, (ii) says that the end
diagrams are switch-free classical link diagrams, (iii) says 
that the classical crossing information is consistent, and (iv) avoids 
forbidden moves. Then condition (v) says that the sequence is maximally 
virtualized.
\end{proof}

Examples of unliftable immersed surfaces are known (see \cite{car} for some 
examples, such as a double cover of Boy's Surface), but it is unknown 
whether any such unliftable surface satisfies the conditions in theorem 
\ref{re-1}.

Each undecorated double point curve corresponds to a virtual $ir$ class, and
sets of undecorated double point curves represent possible realization sets. 
While the surface $M$ has no local extrema in the vertical direction, the 
double-point curves typically do. Each portion of a double point curve 
joining a maximum or minimum on the curve (including endpoints) is the 
``trajectory'' of a virtual crossing, and two portions of a curve meeting 
at a maximum or minimum correspond to $i$ or $r$ partners. Indeed, this 
observation provides another proof of proposition \ref{ir}.

Theorem 4.6 of \cite{car} says that an immersed surface is liftable if and 
only if its decker set can be partitioned into two classes $A$ and $B$, where 
every singular curve is the image of one $A$ curve and one $B$ curve, at 
every branch point an $A$ curve and a $B$ curve meet, and at every triple 
point the preimage in $M$ consists of three intersections of decker curves, 
one involving two $A$ curves, one involving two $B$ curves, and one involving 
one $A$ and one $B$. Hence, To prove conjecture \ref{main} it would suffice 
to show that every immersed surface digram meeting the conditions in the 
theorem \ref{re-1} above admits such a coloring of the virtual double decker 
curves by $A$ and $B$ compatible with the crossing information specified 
for all classical decker curves intersecting the virtual decker curves in 
triple points.

A curve in an immersed surface diagram $M$ satisfying the condition of theorem
\ref{re-1} which misses the double point set and connects the end diagrams 
with no relative extrema in the vertical direction defines a base point fixed 
by the moves in the virtual isotopy sequence; call such a curve a \textit{base
curve}. If $M$ is connected, then the vertical slices which miss triple points
are virtual knot diagrams; hence, if there is a base curve and a choice of 
direction of travel around the diagrams starting at the base point such that 
the preimage of every horizontal slice of $M$ meets the upper decker curve 
of each classical double curve intersecting two virtual curves at a triple 
point before the lower decker curve, then by theorem \ref{vd}, the 
surface is liftable. 

\section{Knots in Surfaces}

In \cite{cks}, virtual links are shown to be equivalent to links in 
thickened surfaces $S\times I$ with a stabilization operation consisting 
of adding or removing handles which miss the link. This corresponds to the 
intuitive concept of virtual knot diagrams as non-planar link diagrams; if we
draw our link $L$ on a surface $S$ and then project $S$ onto $\mathbb{R}^2$,
virtual crossings arise as the result of parts of the link in distinct handles 
or opposite sides of a handle projecting to the same point in the plane.

More specifically, a virtual link diagram $D\subset \mathbb{R}^2$ is the 
image of a link $L\subset S\times I$ under the composition $p=p_2 \circ p_1$
where $p_1:S\times I\to S$ and $p_2: S\to \mathbb{R}^2$. Virtual crossings
and moves appear only in the projection $p_2$ and hence are dependent on
the choice of embedding of $S\times I$ in $\mathbb{R}^3$. For a given 
$L\subset S\times I$, there may be many non-isotopic embeddings of 
$S\times I$ into $\mathbb{R}^3$, involving knotted and linked handles as 
well as Dehn twists. A sequence of Reidemeister and stabilization moves then 
corresponds to a sequence of choices of embedding 
$S\times I\hookrightarrow \mathbb{R}^3 $ and projection $p=p_2 \circ p_1$ 
which yields a virtual isotopy sequence, with the classical crossings
occurring as double points of $p_1$ and the virtual crossings as double points
of $p_2$.

In \cite{ku}, it is shown that every stable equivalence class of links in 
thickened surfaces has a unique irreducible representative; that is, given 
a link $L$ in a thickened surface $S\times I$, any two sequences of 
destabilization moves resulting in irreducible surfaces yield homeomorphic 
$(S \times I, L)$ pairs. In particular, any two sequences of Reidemeister 
moves and stabilization moves on $L\subset S\times I$ which starts and ends 
with genus zero, i.e. classical, diagrams, result in diagrams which differ 
only by Reidemeister moves.

An embedding of a thickened surface containing a link into $\mathbb{R}^3$ is 
a virtual crossing realization for the corresponding virtual link diagram, 
since the choice of embedding for the surface includes a choice of over/under 
for each handle; we may view the virtual crossing realization as simply 
forgetting the distinction between crossings arising from $p_1$ and $p_2$. 
Conversely, stabilization moves in an embedded surface in $\mathbb{R}^3$ can 
result in virtualizing virtualizable crossings. A set of classical crossings 
is sequentially virtualizable if the crossings can be removed via 
stabilization.

\begin{definition}
Let $L\subset S\times I$ be a link in a thickened surface.
An embedding $e:S\times I\hookrightarrow \mathbb{R}^3$ is a 
\textit{lift} of a virtual knot diagram $D$ if $D=p(L)$ where 
$p=p_2 \circ p_1:e(S\times I) \to e(S\times 0)\to \mathbb{R}^2.$ 
\end{definition}

With this definition, we end with another reformulation of conjecture 
\ref{main}, namely:

\begin{theorem}
Conjecture \ref{main} is equivalent to: 
Every maximally virtualized virtual isotopy sequence with switch-free 
classical end diagrams $K_1\to \dots \to K_n$ has a sequence of lifts 
$e_i:S\times I \hookrightarrow \mathbb{R}^3,\ i=1\dots n$ 
such that $e_i$ is a lift of $K_i$ and $e_i(L)$ is ambient isotopic 
to $e_{i+1}(L)$ in $\mathbb{R}^3$ for each $i=1\dots n$ by an isotopy which
fixes the part of the link fixed by the corresponding virtual move.
\end{theorem}

\begin{proof}
If a maximally virtualized virtual isotopy sequence with switch-free 
classical end diagrams admits a virtual crossing realization, this realization
tells us how to choose the lifts for each diagram so that the link $L$ is 
changed only by ambient isotopies in $\mathbb{R}^3$ by specifying handle 
crossing and other embedding information. Conversely, if every such sequence 
is liftable, the specific lifts define a virtual crossing realization, and 
the condition that the links are ambient isotopic guarantees that the 
realization is valid.
\end{proof}

\section{Acknowledgments}

The author is grateful to several people whose conversations were most helpful
during the preparation of this paper, including but not limited to 
R.~A.~Litherland, Louis Kauffman, Scott Carter, Masahico Saito, 
Shin Satoh, who suggested that the author consider the problem in terms of 
2-knots, and the reviewer, whose helpful comments and suggestions have 
improved the paper significantly.

\nonumsection{References}

\end{document}